\documentclass[11pt,fullpage]{article}
\usepackage{amsmath,amssymb,amsfonts,amsthm}
\newcommand{\be}{\begin{equation}}
\newcommand{\ee}{\end{equation}}
\newcommand{\bea}{\begin{eqnarray}}
\newcommand{\eea}{\end{eqnarray}}
\newcommand{\bean}{\begin{eqnarray*}}
\newcommand{\eean}{\end{eqnarray*}}
\newcommand{\ben}{\begin{equation}{nonumber}}
\newcommand{\een}{\end{equation}{nonumber}}

\numberwithin{equation}{section}

\newtheorem{dfn}{Definition}[section]
\newtheorem{thm}[dfn]{Theorem}
\newtheorem{lema}[dfn]{Lemma}
\newtheorem{pro}[dfn]{Proposition}
\newtheorem{coro}[dfn]{Corollary}
\newtheorem{xmpl}[dfn]{Example}
\newtheorem{rmrk}[dfn]{Remark}
\newcommand{\bdfn}{\begin{dfn}}
\newcommand{\bthm}{\begin{thm}}
\newcommand{\blema}{\begin{lema}}
\newcommand{\bpro}{\begin{pro}}
\newcommand{\bcoro}{\begin{coro}}
\newcommand{\bxmpl}{\begin{xmpl}}
\newcommand{\brmrk}{\begin{rmrk}}
\newcommand{\edfn}{\end{dfn}}
\newcommand{\ethm}{\end{thm}}
\newcommand{\elema}{\end{lema}}
\newcommand{\epro}{\end{pro}}
\newcommand{\ecoro}{\end{coro}}
\newcommand{\exmpl}{\end{xmpl}}
\newcommand{\ermrk}{\end{rmrk}}

\newcommand{\B}{\mathcal B}
\newcommand{\A}{\mathcal A}
\newcommand{\E}{\mathcal E}

\newcommand{\K}{\mathcal K}

\newcommand{\raro}{\rightarrow}

\begin{document}
\begin{center}
 {\large {\bf Quantum random walks and vanishing of the second
 Hochschild cohomology}}\\
 {\large  Debashish Goswami  }\\
{Stat-Math Unit, Indian Statistical Institute, Kolkata,\\ 203, B.T.
Road, Kolkata-108,
 India.\\ E-mail :
goswamid@isical.ac.in}\\
 {\large  Lingaraj Sahu {\footnote {Author
acknowledge the support
of National Board of Higher Mathematics, DAE, India .}}}\\
{Stat-Math Unit, Indian Statistical Institute, Bangalore Centre,\\
$8^{th}$
 Mile, Mysore Road,
 Bangalore-59,
 India.
 \\E-mail : lingaraj@isibang.ac.in}
\end{center}
\begin{abstract}
Given a conditionally completely positive map $\mathcal L$ on a
unital  $\ast$-algebra $\A$, we find an interesting connection
between the second Hochschild cohomology of  $\A$ with coefficients
in the bimodule $E_{\mathcal L}=\B^a(\A \oplus M)$ of adjointable
maps, where $M$ is the GNS bimodule of $\mathcal L$,  and the
possibility of constructing a quantum random walk (in the sense of
\cite{AP,LP,L,KBS}) corresponding to $\mathcal L$.
\end{abstract}
 \section{Introduction}
Quantum dynamical semigroups (QDS for short), which are $C_0$-
semigroups of completely positive, contractive maps on $C^*$ or von
Neumann algebras (with appropriate continuity assumptions), are
interesting and important objects of study both from physical as
well as mathematical viewpoints. A very useful tool for
understanding such semigroups is Evans-Hudson dilation (EH dilation
for short). By an E-H dilation of a QDS $(T_t)_{t \geq 0}$ on a von
Neumann algebra $\A \subseteq {\mathcal B}(\mathbf{h})$, we mean a
family $j_t$ of normal $\ast$-homomorphism from $\A$ into $\A
\otimes {\mathcal B}(\mathbf{h} \otimes \Gamma(L^2(R_+,\mathbf
k)))$, where $\mathbf k$ is a Hilbert space, and $j_t$ satisfies a
quantum stochastic differential equation of the form
$$ dj_t(x)=j_t(\theta^\alpha_\beta(x)),~~j_0(x)=x \otimes I,$$ for
$x$ belonging to a suitable dense $\ast$-subalgebra on which a
family of linear maps $\theta^\alpha_\beta$ are defined, and
$\theta^0_0$ coincides with the generator of $T_t$. For more details
of this concept, we refer the reader to the books \cite{KRP,GS} and
references therein. While there is a complete theory of such
dilations for semigroups with norm-bounded generator (i.e. uniformly
continuous semigroups), there is hardly any hope for a general
theory for an arbitrary QDS. Nevertheless, there have been several
attempts to construct EH dilation for different classes of QDS with
unbounded generator. Moreover, there are more than one constructions
of the family $j_t$ for a QDS with bounded generator. In addition to
the traditional approach by iteration, there is a very interesting
construction (see \cite{L,S,Bel2}) of EH dilation as a strong limit
of a sequence of homomorphism which can be thought of as `quantum
random walk'. It should be mentioned that for building a
satisfactory general theory of EH dilation covering a reasonably
large class of QDS with unbounded generator, it is absolutely
crucial to deeply look into all the different approaches available
in the bounded generator case, and to see whether some of them, or a
suitable combination of them, can be generalized to cover QDS with
unbounded generator. Indeed, the approach through quantum random
walk seems to have a great promise in this context. However, there
are two issues involved in this approach : first, to construct a
quantum random walk for a given QDS (possibly with unbounded
generator), and then to see whether it converges strongly. In the
present article we study some algebraic conditions for the
possibility of constructing a quantum random walk in the general
situation. We work in a purely algebraic setting, and are able to
discover a very interesting connection between  the algebraic
relations satisfied by components of a quantum random walk (if it
exists) and the second Hochschild cohomology of the algebra with
coefficient in a module naturally associated with the CCP generator
of the given QDS. We leave the study of analytic aspects of our
results for later work. It may be remarked here that the first and
second Hochschild cohomolgies did appear in several other works on
QDS and quantum probabilistic dilation, for example the celebrated
work of Christensen and Evans (\cite{CE}), and also in work of
Hudson (\cite{Hud}). However, none of those works are concerned with
the quantum random walks and do not have any overlap with the
results obtained in the present article.

\section{Notations and Preliminaries}

{\bf Quantum Random Walk}

Let ${\mathcal K} =L^2(\mathbb{R}_+,{\mathbf k})$  where ${\mathbf
k}$ is a Hilbert space and let $\Gamma$ be the  symmetric Fock space
$\Gamma({\mathcal K})$ over $\mathcal K.$ For any partition
$S\equiv(0=t_0<t_1<t_2\cdots)$ of $\mathbb{R}_+,~\mathcal
K=\oplus_{n\ge 1} {\mathcal K}_n,$  where  ${\mathcal K}_n$ is the
range  of projection $1_{(t_{n-1},t_n]}$ and  the Fock space
$\Gamma$ can be viewed as the infinite tensor product $\otimes_{n\ge
1}\Gamma_n$ of symmetric Fock spaces $\{\Gamma_n=\Gamma({\mathcal
K}_n)\}_{n\ge 1}$ with respect to the stabilizing sequence
$\Omega=\{ \Omega_n: n\ge 1\},$ where
$\Omega_n=\Omega_{(t_{n-1},t_n]}$ is the vacuum vector in
$\Gamma_n.$ Let denote the interval $(t_{n-1},t_n]$ by $[n]$  and
the orthogonal projection  of $\Gamma_n$ onto the $m$-particle space
by $P_m [n].$

For $n\ge 1,$ consider the subspace $\hat{\mathbf k}_n={\mathbb
{C}}\ \Omega_n\oplus \mathbf k_n$ of $\Gamma,$ where $\mathbf
k_n=\{1_{[n]}\phi: \phi\in \mathbf k \}.$ The spaces $ \hat{\mathbf
k}_n,$ are  isomorphic with  $ \hat{\mathbf k}:={\mathbb {C}}\oplus
{\mathbf k} .$ \bdfn The toy Fock space associated with the
partition $S$ of ${\mathbb R}_+$  is defined to be the subspace
$\Gamma(S):=\otimes_{n\ge 1} \hat{\mathbf k}_n$ with respect to the
stabilizing sequence  $( \Omega_n)_{ n\ge 1}.$ \edfn

Let $P(S)$ be the orthogonal projection of $\Gamma$ onto the toy
Fock space $\Gamma(S).$ Now onwards let us consider toy Fock space
$\Gamma(S_h)$ associated with  regular partition  $S_h\equiv
(0,h,\cdots)$ for some   $h>0$ and denote the orthogonal projection
by $P_h.$  Denoting the restriction of orthogonal projection $P_h$
 to $\Gamma_n$ by $P_h[n], P_h=\otimes_{n\ge 1} P_h[n].$
Now we define basic operators associated with toy Fock space
$\Gamma(S_h)$ using the fundamental processes in coordinate-free
language of quantum stochastic calculus, developed in \cite{GS1}.
For $S\in \B(\mathbf h), R\in \B(\mathbf h,\mathbf h\otimes {\mathbf
k}),$ $Q\in \B(\mathbf h\otimes {\mathbf k},\mathbf h)$ and $T\in
\B(\mathbf h\otimes {\mathbf k})$ let us define four basic operators
on  $\Gamma$ as follows, for $n\ge 1,$
\begin{equation}
\begin{split}
&N_{S}^{00}[n]=S P_0[n] ,\\
&N_{Q}^{01}[n]= \frac{a_Q[n]}{\sqrt{h}}P_1[n],\\
&N_{R}^{10}[n]=P_1[n] \frac{a^\dag_R[n]}{\sqrt{h}}, \\
&N_{T}^{11}[n]=P_1[n](\lambda_T[n])P_1[n]P_h[n].
\end{split}
\end{equation}
Here all these operators  act nontrivially only  on $\Gamma_n.$
 For definition of coordinate-free fundamental
processes $\Lambda$'s we refer to \cite{GS1}. Here, we note that in
the notation of \cite{GS1}, the annihilation process $ a_Q[n]$
appear above is $a_{Q^*}[n].$ All these maps  $  \B(\mathbf h)\ni
S\mapsto N_S^{00}[n],$ $  \B(\mathbf h \otimes {\mathbf k},\mathbf
h)\ni Q\mapsto N_Q^{01}[n],$
 $ \B(\mathbf
h,\mathbf h \otimes {\mathbf k})\ni
 R\mapsto N_R^{10}[n] $
and $ \B(\mathbf h \otimes {\mathbf k})\ni T\mapsto N_T^{11}[n]$ are
linear. It is clear that these   operators $N$'s are bounded and
leave the subspace $\Gamma(S_h)$  invariant.
 It can be shown that (for detain see
\cite{S}):\\
\begin{itemize} \label{N-homo}
\item $(N_X^{\mu \nu}[n])^*=N_{X^*}^{\nu \mu}[n],~\forall \mu,\nu\in\{0,1\}$

\item $N_{S}^{00}[n]+N_{S\otimes 1_{\mathbf k}}^{11}[n]=S\otimes P_h[n]$

\item $N_X^{\mu \nu}[n]N_{Y}^{\eta \xi}[n]=\delta_\nu^\eta N_{XY}^{\mu
\xi}[n],$  where $\delta_\nu^\eta$ is Dirac delta function of $\eta$
and $\nu.$

\end{itemize}

 Let $\A$ be a unital $*$-subalgebra of $\B(\mathbf h).$
Suppose we are given with a family of $*$-homomorphisms
$\{\beta(h)\}_{h > 0}$ from $\A$ to $\A\otimes\B(\hat{\mathbf k}).$
It can be written that $\beta(h)=\left(
\begin{array}{cc} \beta_{00}(h)
  & \beta_{01}(h)\\
\beta_{10}(h) &  \beta_{11}(h)
\end{array} \right),$  where the components $\beta_{00}(h)\in \B(\A),~
\beta_{11}(h)\in \B(\A,
 \A\otimes
\B({\mathbf k}))$ and $\beta_{10}(h) \in \B(\A,\A\otimes {\mathbf
k})$ such that
 \bean
  &&\beta_{\mu\nu}(h)(x^*)=(\beta_{\nu\mu}(h)(x))^*, \\
&&\beta_{\mu\nu}(h)(xy)=\sum_{\eta=0}^1\beta_{\mu\eta}(h)(x)\beta_{\eta
\nu}(h)(y). \eean
 Let us define a family of  maps $\{{\mathcal
P}_t^{(h)}:\A \otimes \E(\K)\raro \A\otimes \Gamma\}_{t\ge 0}$ as
follows. First subdivide the interval $[0,t]$ into $[k]\equiv \left(
(k-1)h,kh\right],~1\le k\le n$ so that $t\in \left(
(n-1)h,nh\right]$  and set for $x\in \A,~f\in \K$
\begin{equation}
\left.
\begin{array}{l}
{\mathcal P}_{0}^{(h)}(x \textbf{e}(f))=x \textbf{e}(f)\\ \\
{\mathcal P}_{kh}^{(h)}(x \textbf{e}(f))=\sum_{\mu,\nu=0}^1
{\mathcal P}_{(k-1)h}^{(h)}  N_{\beta_{\mu\nu}(h,x)}^{\mu\nu}[k]
\textbf{e}(f)
\end{array}
\right \}
\end{equation}
and
${\mathcal P}_{t}^{(h)}={\mathcal P}_{nh}^{(h)}.$\\
Setting $p_t^{(h)}(x)u \textbf{e}(f):={\mathcal P}_t^{(h)}(x
\textbf{e}(f))u, \forall u\in\mathbf h,$ by the properties of the
family $\{\beta_{\mu\nu}(h)\}$ and $\{N^{\mu\nu}[k]\},$ $p_t^{(h)}$
are $*$-homomorphism from $\A$ into $\A\otimes \B(\Gamma).$ \bdfn
This family of $*$-homomorphisms $\{p_t^{(h)}:t\ge 0\}$ is called  a
quantum random walk (QRW) associated with $\beta(h).$ \edfn
\noindent{\bf
Hochschild cohomology} \\
Let us recall the definition of the Hochschild cohomology
$H^n(\A,N)$ for $\A$ with coefficients in an $\A$-$\A$ bimodule $N$
(for detail we refer to \cite{SiSm}). It is the cohomology of the
cochain complex $(C^n\equiv C^n(\A,N), b),~n \geq 0$, where $C^0=N$,
and for $n \geq 1,$~ $C^n$ consists of
 all multi-$\mathbb C $-linear maps $f : \A \times \cdots \A~ (~n~{\rm copies})
 \rightarrow N$, with
 the coboundary map $b$ given by \\
 $ bf(a_0,a_1,\cdots,a_n)
 :=a_0f(a_1,\cdots,a_n)
 $
 $$+\sum_{i=0}^{n-1} (-1)^{i+1} f(a_0,\cdots,a_{i-1},a_ia_{i+1},
 \cdots,a_n) +(-1)^{n+1}a_nf(a_0,\cdots,a_n).$$
 Let us introduce one more notation. Let ${\mathcal C}[[ t]]$
 denote the ring of formal power series
  in one indeterminate $t$ with coefficients in a ring ${\mathcal C}$.
  If $\mathcal C$ is a $\ast$-algebra,
  so is ${\mathcal C}[[t]]$.

 \section{Main results}
Let $\A$ be a unital $*$-subalgebra of  $ \B(\mathbf{h})$  and
${\mathcal L}$ be a  conditionally  completely positive (CCP)  map
from  $\A$ into itself, satisfying ${\mathcal L}(1)=0$.  Then there
exist  a canonical (unique upto isomorphism) pre-Hilbert $\A$-$\A$
bimodule  $M$, with the left action denoted by $\pi$ (can also be
viewed as a $\ast$-representation  $\pi$  of  the $*$-algebra $\A$
into the  algebra $\B^a(M)$ of adjointable maps  on $M$ ), and a
bimodule-derivation $\delta : \A \rightarrow M$, such that $M$
coincides with the right $\A$-linear span of $\delta(\A)$. Note that
we can identify an element $\xi \in M$ with the rank-one map $\xi^*
\equiv < \xi| : M \rightarrow \A$ given by $M \ni \eta \mapsto
<\xi,\eta> \in \A$, where $<\cdot, \cdot>$ denotes the $\A$-valued
inner product on $M$. We have
\[\mathcal L(xy)- x\mathcal L(y)-\mathcal L(x)y=
\delta^\dag(x)\delta(y), ~~\forall
 x,y \in  \A,\] where $\psi^\dag$ for a linear map $\psi$ on $\A$
 is defined as $\psi^\dag(x):= (\psi(x^*))^*.$

 When $\A$ is  a von Neumann algebra and $\mathcal L$ is norm-bounded then
 one can imbed $M$ in a Hilbert von
  Neumann module of the form
  ${\mathcal A} \otimes \mathbf k\subseteq(\B(\mathbf h,\mathbf h\otimes \mathbf k))$
  for some Hilbert space $\mathbf k$, and show that $\delta, \delta^\dagger,
  \pi$ are all bounded maps. Furthermore, using the explicit structure of $\mathcal L$
  as obtained from the Christensen-Evans
  Theorem (ref. \cite{CE}) one can construct (see \cite{Bel2,S}) a family of
   $*$-homomorphism $\{\beta(h):\A \raro \A\otimes \B(\hat{\mathbf k}) ~:~h > 0\}$
 such that
 $\beta(h)=\left( \begin{array}{cc} \beta_{00}(h)
  & \beta_{01}(h)\\
\beta_{10}(h) &  \beta_{11}(h)
\end{array} \right),$ where
\begin{itemize}
\item $\beta_{00}(h)= \sum_{n\ge 0}  h^n \theta_{00}^{(n)}$ with $\theta_{00}^{(0)}(x)=x,
 \theta_{00}^{(1)}(x)=\theta_{00}(x)$
\item $\beta_{10}(h)= \sum_{n\ge 1}  h^{\frac{2n-1}{2}} \theta_{10}^{(n)}$
with $\theta_{10}^{(1)}(x)=\delta(x),$
\item $\beta_{10}(h)= \sum_{n\ge 1}  h^{\frac{2n-1}{2}} \theta_{01}^{(n)}$
with $\theta_{01}^{(1)}(x)=\delta^\dag(x),$
\item $\beta_{11}(h)= \sum_{n\ge 1}  h^{n-1} \theta_{11}^{(n)}$ with
$\theta_{11}^{(1)} (x) = \pi(x).$
 \end{itemize}
 Using this, an EH flow  for the QDS generated by
  $\mathcal L$ can be constructed (see \cite{Bel2,S}) as strong limit of quantum random
  walks discussed in the previous section.

However, in this paper we concentrate on the purely algebraic
aspect of such construction only and make the interesting
observation that this is related intimately to the vanishing of
second Hochschild cohomology of $\A$.

 Now, for a purely algebraic treatment, let us fix a  $\ast$-algebra $\A$, CCP map
 $\mathcal L$ as in the beginning, and the bimodule $M$ and the derivation $\delta$
 as mentioned before (not assumed to be bounded in any sense).
 Let us also consider the pre-Hilbert $\A$-$\A$ bimodule $E_{\mathcal L}
 :=\B^a( \hat{M}) \equiv \B^a(\A \oplus M)$, with
the bimodule actions given by   $x.R=\tilde{\pi}(x)R $  and $R.x= R
\tilde{\pi}(x)$ where $\tilde{\pi}(x) =x\oplus  \pi(x) $. Let us
denote by $M^\ast$ the submodule of $E_{\mathcal L}$ consisting of
$\xi^*, ~\xi \in M$. It is clear that $M$, $M^\ast$, $\A$ and
$\B^a(M)$ are canonically imbedded as complemented submodules of
$E_{\mathcal L}$ and in fact, $E_{\mathcal L}$ is the direct sum of
these four submodules. Any element  $X $ of $ E_{\mathcal L}$ can be
written as a $2 \times 2$ matrix
form $$ \left( \begin{array}{cc}X_{11} &  X_{12}\\ X_{21} &  X_{22} \\
\end{array} \right),$$
where $X_{11} \in \A$, $X_{12} \in M^*$, $X_{21} \in M$ and $X_{22}
\in \B^a(M)$.

\bthm \label{main} If $H^2(\A,E_{\mathcal L})=0$ then  there exists
a $*$-homomorphism $\beta:\A \raro E_{\mathcal L} [[ t]]$ such that
$\beta(t)=\left( \begin{array}{cc} \beta_{00}(h)
& \beta_{01}(h)\\
\beta_{10}(h) &  \beta_{11}(h)
\end{array} \right),$ where $h=t^2$ and
\begin{itemize}
\item $\beta_{00}(h)= \sum_{n\ge 0}  h^n \theta_{00}^{(n)}$ with $\theta_{00}^{(0)}(x)=x,
\theta_{00}^{(1)}(x)=\theta_{00}(x)$
\item $\beta_{10}(h)= \sum_{n\ge 1}  h^{\frac{2n-1}{2}} \theta_{10}^{(n)}$
with $\theta_{10}^{(1)}(x)=\delta(x),$
\item $\beta_{10}(h)= \sum_{n\ge 1}  h^{\frac{2n-1}{2}} \theta_{01}^{(n)}$
with $\theta_{01}^{(1)}(x)=\delta^\dag(x),$
\item $\beta_{11}(h)= \sum_{n\ge 1}  h^{n-1} \theta_{11}^{(n)}$ with
$\theta_{11}^{(1)} (x) = \pi(x).$
\end{itemize}

\ethm
\begin{proof}
First of all we note that $H^2(\A, N)$ =0 for any complemented
submodule $N$ of $E_{\mathcal L}$, for example, for $N=M, M^\ast,
\A, \B^a(M)$. Moreover, we shall view any map from some module to
any such submodule $N$ of $E_{\mathcal L}$ as a map into
$E_{\mathcal L}$. Also, it is easy to verify that  the
$*$-homomorphic property of $\beta$ is equivalent to
 \bea
 \label{*homo}
  &\beta_{\mu\nu}(h)(x^*)=(\beta_{\nu\mu}(h)(x))^*, \\
&\beta_{\mu\nu}(h)(xy)=\sum_{\eta=0}^1\beta_{\mu\eta}(h)(x)\beta_{\eta
\nu}(h)(y). {\nonumber} \eea

To prove existence  and $*$-homomorphic properties of $\beta,$ by
induction, we shall show the existence  of maps
$\theta_{\mu\nu}^{(n)}\in \mathcal C^1(\A, E_{\mathcal L} )$
satisfying
\begin{equation}
\begin{split} \label{algebraic-relation}
&\theta_{11}^{(n)}(xy)=\sum_{k=1}^{n-1} \theta_{10}^{(k)}(x)
\theta_{01}^{(n-k)}(y) + \sum_{k=1}^{n} \theta_{11}^{(k)}(x)
\theta_{11}^{(n-k+1)}(y)\\
&\theta_{10}^{(n)}(xy)= \sum_{k=1}^{n} \theta_{10}^{(k)}(x)
\theta_{00}^{(n-k)}(y)+ \sum_{k=1}^{n} \theta_{11}^{(k)}(x)
\theta_{10}^{(n-k+1)}(y)\\
&\theta_{01}^{(n)}(xy)= \sum_{k=0}^{n-1} \theta_{00}^{(k)}(x)
\theta_{01}^{(n-k)}(y) +\sum_{k=1}^{n}
 \theta_{01}^{(k)}(x)\theta_{11}^{(n-k+1)}(y)\\
&\theta_{00}^{(n)}(xy) =\sum_{k=0}^{n}\theta_{00}^{(k)}(x)
\theta_{00}^{(n-k)}(y) +\sum_{k=1}^{n}
 \theta_{01}^{(k)}(x) \theta_{10}^{(n-k+1)}(y) \\
&\theta_{\mu \nu}^{(n)}(x^*)=(\theta_{\nu \mu}^{(n)}(x))^*.
\end{split}
\end{equation}

First, let us consider the following elements of $\mathcal C^2(\A,
E_{\mathcal L} )$  and $\mathcal C^1(\A, E_{\mathcal L} )$

 \begin{itemize}
\item $\phi_{11}^{(2)}(x,y):= \theta_{10}^{(1)}(x)
\theta_{01}^{(1)}(y).$ \\
As $\partial \theta_{10}^{(1)}(x,y)=0,~\partial
\theta_{01}^{(1)}(x,y)=0$ we have
$\partial\phi_{11}^{(2)}(x,y,z)=0.$ Now since $H^2(\A,E_{\mathcal
L})=0,$ there exists  a map, say $\theta_{11}^{(2)}\in \mathcal
C^1(\A,E_{\mathcal L})$ such that $\partial
\theta_{11}^{(2)}=\phi_{11}^{(2)}.$ Since we have
${\theta_{01}^{(1)}}^\dagger=\theta_{10}^{(1)}$, it is easy to see
that $(\phi_{11}^{(2)}(y^*,x^*))^*=\phi_{11}^{(2)}(x,y)$, so
$\partial {\theta_{11}^{(2)}}^\dagger=\partial \theta_{11}^{(2)}$.
 Thus,  taking $\gamma= \frac{1}{2} (\theta_{11}^{(2)}+{\theta_{11}^{(2)}}^\dagger)$,
 we have
 $\partial \gamma=\phi_{11}^{(2)}$ and $\gamma^\dagger=\gamma$.
  By replacing $\theta_{11}^{(2)}$ by $\gamma$, we can assume without loss of
  generality that $\theta_{11}^{(2)}(x^*)^*=\theta_{11}^{(2)}(x)$.

\item $\phi_{10}^{(2)}(x,y):= \theta_{10}^{(1)}(x) \theta_{00}^{(1)}(y)+
\theta_{11}^{(2)}(x) \theta_{10}^{(1)}(y).$
 Now \bean
&&\partial \phi_{10}^{(2)}(x,y,z)\\
&&= \theta_{10}^{(1)}(x) \partial \theta_{00}^{(1)}(y,z)-\partial
\theta_{10}^{(1)}(x,y) \theta_{00}^{(1)}(z)+ \theta_{11}^{(2)}(x)
\partial \theta_{10}^{(1)}(y,z)-\partial \theta_{11}^{(2)}(x,y)
\theta_{10}^{(1)}(z)\\
&&=\theta_{10}^{(1)}(x) \theta_{01}^{(1)}(y)
\theta_{10}^{(1)}(z)-0+0-\theta_{10}^{(1)}(x) \theta_{01}^{(1)}(y)
\theta_{10}^{(1)}(z)=0 .\eean Since $H^2(\A,E_{\mathcal L})=0,$
there exists a map, say $\theta_{10}^{(2)}\in \mathcal
C^1(\A,E_{\mathcal L})$ such that $\partial
\theta_{10}^{(2)}=\phi_{10}^{(2)}.$

Now define $ \theta_{01}^{(2)}(x):= (\theta_{10}^{(2)}(x^*))^*.$
 Then \bean &&\phi_{01}^{(2)}(x,y):=\partial \theta_{01}^{(2)}(x,y)=
\theta_{01}^{(2)}(xy)-x\theta_{01}^{(2)}(y)-\theta_{01}^{(2)}(x) \pi
(y) \\
&&=\{\theta_{10}^{(2)}(y^* x^*)-\theta_{10}^{(2)}(y^*) x^*- \pi
(y^*)\theta_{10}^{(2)}(x^*) \pi
(y)\}^*=\{\partial \theta_{10}^{(2)}(y^*, x^*)\}^*\\
&&=\{ \theta_{10}^{(1)}(y^*) \theta_{00}^{(1)}(x^*)+
\theta_{11}^{(2)}(y^*) \theta_{10}^{(1)}(x^*)\}^*\\
 &&=\theta_{00}^{(1)}(x) \theta_{01}^{(1)}(y) +
 \theta_{01}^{(1)}(x)\theta_{11}^{(2)}(y).
 \eean

\item $\phi_{00}^{(2)}(x,y):=\theta_{00}^{(1)}(x) \theta_{00}^{(1)}(y)
+ \theta_{01}^{(1)}(x) \theta_{10}^{(2)}(y)+ \theta_{01}^{(2)}(x)
\theta_{10}^{(1)}(y) .$
Now
 \bean
&&\partial \phi_{00}^{(2)}(x,y,z)\\
&&=\theta_{00}^{(1)}(x)\partial \theta_{00}^{(1)}(y,z)
-\partial\theta_{00}^{(1)}(x,y) \theta_{00}^{(1)}(z) \\
&&~~~ + \theta_{01}^{(1)}(x)
\partial\theta_{10}^{(2)}(y,z)-\partial\theta_{01}^{(1)}(x,y)
\theta_{10}^{(2)}(z)\\
&&~~+ \theta_{01}^{(2)}(x)
\partial\theta_{10}^{(1)}(y,z) -\partial\theta_{01}^{(2)}(x,y)
\theta_{10}^{(1)}(z) \\
&&= \theta_{00}^{(1)}(x)\theta_{01}^{(1)}(y) \theta_{10}^{(1)}(z)
-\theta_{01}^{(1)}(x)
\theta_{10}^{(1)}(y) \theta_{00}^{(1)}(z) \\
&&~~~ + \theta_{01}^{(1)}(x) \{\theta_{10}^{(1)}(y)
\theta_{00}^{(1)}(z)+
\theta_{11}^{(2)}(y) \theta_{10}^{(1)}(z)\}-0\\
&&~~+0-\{ \theta_{00}^{(1)}(x) \theta_{01}^{(1)}(y) +
 \theta_{01}^{(1)}(x)\theta_{11}^{(2)}(y)\} \theta_{10}^{(1)}(z)\\
 &&=0. \eean
 Since $H^2(\A, E_{\mathcal L} )=0,$  there exists a
map, say $\theta_{00}^{(2)}\in \mathcal C^1(\A, E_{\mathcal L} )$
such that
$\partial \theta_{00}^{(2)}=\phi_{00}^{(2)}.$\\
As seen before, it can be arranged, by replacing
$\theta_{00}^{(2)}$ by
$\frac{1}{2}(\theta_{00}^{(2)}+{\theta_{00}^{(2)}}^\dagger)$ if
necessary,  that $\theta_{00}^{(2)}(x^*)=(\theta_{00}^{(2)}(x))^*$
\end{itemize}

Now we  prove by induction that  there exists a family of  maps
$\{\theta_{\mu\nu}^{(n)}\in \mathcal C^1(\A, E_{\mathcal L}
):\mu,\nu\in \{0,1\}, n\ge 1 \}$   such that

 \begin{enumerate}
 \item $\partial \theta_{11}^{(n)}(x,y)=\theta_{11}^{(n)}(xy)-\pi(x)\theta_{11}^{(n)}(y)
-\theta_{11}^{(n)}(x) \pi(y)\\
= \sum_{k=1}^{n-1} \theta_{10}^{(k)}(x) \theta_{01}^{(n-k)}(y) +
\sum_{k=2}^{n-1} \theta_{11}^{(k)}(x) \theta_{11}^{(n-k+1)}(y)$\\
$\theta_{11}^{(n)}(x^*)=(\theta_{11}^{(n)}(x))^*$

\item $\partial \theta_{10}^{(n)}(x,y)=\theta_{10}^{(n)}(xy)-\pi(x)\theta_{10}^{(n)}(y)
-\theta_{10}^{(n)}(x) y\\
=\sum_{k=1}^{n-1} \theta_{10}^{(k)}(x) \theta_{00}^{(n-k)}(y)+
\sum_{k=2}^{n} \theta_{11}^{(k)}(x) \theta_{10}^{(n-k+1)}(y)$

\item  $\theta_{01}^{(n)}(x)=(\theta_{01}^{(n)}(x^*))^* ,\\
\partial \theta_{01}^{(n)}(x,y)=\theta_{01}^{(n)}(xy)-x\theta_{01}^{(n)}(y)
-\theta_{01}^{(n)}(x) \pi(y)\\= \sum_{k=1}^{n-1}
\theta_{00}^{(k)}(x) \theta_{01}^{(n-k)}(y) +\sum_{k=1}^{n-1}
 \theta_{01}^{(k)}(x)\theta_{11}^{(n-k+1)}(y)$

\item $\partial \theta_{00}^{(n)}(x,y)
=\theta_{00}^{(n)}(xy)-x\theta_{00}^{(n)}(y)-\theta_{00}^{(n)}(x)y
\\=\sum_{k=1}^{n-1}\theta_{00}^{(k)}(x) \theta_{00}^{(n-k)}(y)
+\sum_{k=1}^{n}
 \theta_{01}^{(k)}(x) \theta_{10}^{(n-k+1)}(y) .$\\
  $\theta_{00}^{(n)}(x^*)=(\theta_{00}^{(n)}(x))^*.$

\end{enumerate}
 Let us assume that for some $m\ge 2,$ there exist a
family of maps  $\{\theta_{\mu\nu}^{(n)}\in \mathcal C^1(\A,
E_{\mathcal L} ):\mu,\nu\in \{0,1\}, n< m \}$ satisfying   above
relations.\\
 Consider the map

 $\phi_{11}^{(m)}(x,y)=\theta_{11}^{(m)}(xy)-\pi(x)\theta_{11}^{(m)}(y)
-\theta_{11}^{(m)}(x) \pi(y)\\
= \sum_{k=1}^{m-1} \theta_{10}^{(k)}(x) \theta_{01}^{(m-k)}(y) +
\sum_{k=2}^{m-1} \theta_{11}^{(k)}(x) \theta_{11}^{(m-k+1)}(y).$\\

Then we have
\bean
&&\partial \phi_{11}^{(m)}(x,y,z)\\
&&= \sum_{k=1}^{m-1} \{ \theta_{10}^{(k)}(x) \partial
\theta_{01}^{(m-k)}(y,z)-\partial \theta_{10}^{(k)}(x,y)
\theta_{01}^{(m-k)}(z)\}\\
 &&~~~+ \sum_{k=2}^{m-1} \{\theta_{11}^{(k)}(x) \partial
\theta_{11}^{(m-k+1)}(y,z)- \partial \theta_{11}^{(k)}(x,y)
\theta_{11}^{(m-k+1)}(z)\}\\
&&= \sum_{k=1}^{m-1}  \theta_{10}^{(k)}(x) \{ \sum_{l=1}^{m-k-1}
\theta_{00}^{(l)}(y) \theta_{01}^{(m-k-l)}(z) +\sum_{l=1}^{m-k-1}
 \theta_{01}^{(l)}(y)\theta_{11}^{(m-k-l+1)}(z)\}\\
&&~~~-\sum_{k=1}^{m-1}\{\sum_{l=1}^{k-1} \theta_{10}^{(l)}(x)
\theta_{00}^{(k-l)}(y)+ \sum_{l=2}^{k} \theta_{11}^{(l)}(x)
\theta_{10}^{(k-l+1)}(y)\}
\theta_{01}^{(m-k)}(z)\\
&&~~~+ \sum_{k=2}^{m-1} \theta_{11}^{(k)}(x) \{\sum_{l=1}^{m-k}
\theta_{10}^{(l)}(y) \theta_{01}^{(m-k-l+1)}(z) +
\sum_{l=2}^{m-k} \theta_{11}^{(l)}(y) \theta_{11}^{(m-k-l+2)}(z)\}\\
&&~~~-\sum_{k=2}^{m-1} \{\sum_{l=1}^{k-1} \theta_{10}^{(l)}(x)
\theta_{01}^{(k-l)}(y) + \sum_{l=2}^{k-1} \theta_{11}^{(l)}(x)
\theta_{11}^{(k-l+1)}(y)\} \theta_{11}^{(m-k+1)}(z) \\
&&=0\eean

Since $H^2(\A, E_{\mathcal L} )=0,$  there exists  a map, say
$\theta_{11}^{(m)}\in \mathcal C^1(\A, E_{\mathcal L} )$ such that
$\partial \theta_{11}^{(m)}=\phi_{11}^{(m)}.$ Moreover,
 it is easily seen that $\partial{ \theta_{11}^{(m)}}^\dagger=\partial \theta_{11}^{(m)}$,
 and so without loss of
  generality we can assume that $\theta_{11}^{(m)}(x^*)=(\theta_{11}^{(m)}(x))^*.$
  Proceeding
similarly  it can be shown the existence of  maps
$\theta_{10}^{(m)}(x),\theta_{01}^{(m)}(x)$   and
$\theta_{00}^{(m)}(x)$ with required relations. \\

From this the
algebraic relations (\ref{algebraic-relation}) follow. Now it is
easy to get (\ref{*homo}), which completes the proof.
  \end{proof}
It is interesting to investigate whether the converse of the above
result also holds; i.e. whether vanishing of $H^2(\A,E)$ is
necessary for the existence of a `quantum random walk' in the formal
algebraic sense as in the above theorem. If the converse to Theorem
\ref{main}  holds, then it will give a `quantum probabilistic'
interpretation
 of $H^2(\A, E)$ as the obstruction to construction of a quantum
 random walk. However, in order to meaningfully apply Theorem \ref{main}
 to the theory of EH dilation, one must obtain an appropriate
 analytic version of it, giving conditions for the formal power
 series in the statement of Theorem \ref{main} to converge. We hope to
 take up these questions in a future work.


\begin{thebibliography}{cgs}

\bibitem{AFL}
Accardi, Luigi; Frigerio, Alberto; Lewis, John T. : Quantum
stochastic processes. Publ. Res. Inst. Math. Sci. {\bf 18}, no. 1,
97--133 (1982).


\bibitem{AP}
Attal, Stéphane; Pautrat, Yan : From repeated to continuous quantum
interactions. Ann. Henri Poincaré {\bf 7} , no. 1, 59--104 (2006).


\bibitem {Bel1}  Belton, Alexander C.R.:  Approximation via toy Fock space –
the vacuum-adapted
 viewpoint, submitted to the Proceedings of the 27th Conference on
 Quantum Probability and applications (Notinghamm 2006).

\bibitem {Bel2} Belton, Alexander C.R.:  Random-walk approximation to vacuum
cocycles, preprint available at
 http://lanl.arxiv.org/math.OA/0702700




\bibitem{CE}  Christensen,   E. and  Evans, D. E. : Cohomology of
operator algebras and quantum dynamical semigroups, J. London Math.
Soc. {\bf 20},  358-368  (1979).



\bibitem {FSk}  Franz U. , Skalski, A.:  Approximation of
quantum L\'{e}vy processes by quantum  random walks,  preprint
available at http://lanl.arxiv.org/math.FA/0703339


\bibitem{GS} Goswami, D. and Sinha, K. B. :
``Quantum Stochastic Processes and Non-commutative Geometry"
Cambridge Tracts in Mathematics, Cambridge University Press, 2007.

\bibitem{GS1} Goswami, D. and Sinha, K.\ B. : \  Hilbert modules and stochastic
dilation of a quantum dynamical semigroup on a von Neumann algebra,
Comm. Math. Phys. \ {\bf 205} no. 2,  377--403  (1999).




\bibitem {GPS} Goswami, D., Pal, A.\ K.\ and Sinha, K.\ B. :   Stochastic
dilation of a quantum dynamical semigroup on a separable unital
$C\sp *$-algebra . Infin. Dimens. Anal. Quantum Probab. Relat. Top.
\ {\bf 3}, no. 1, 177--184 ( 2000).




\bibitem {Hud} Hudson,R.L. : Quantum diffusions and cohomology of algebras.
Proceedings of the
1st World Congress of the Bernoulli Society, Vol. 1 (Tashkent,
1986), 479--483, VNU Sci. Press, Utrecht, 1987.




\bibitem {LP} Lindsay, J. M.; Parthasarathy, K. R. :  The passage from random walk to
diffusion in quantum probability. II. Sankhy\=a Ser. A  {\bf 50} ,
no. 2, 151--170 (1988).


 \bibitem {KRP} Parthasarathy, K.\ R. :\  `` An Introduction to Quantum Stochastic
Calculus'', Monographs in Mathematics, {\bf 85}, Birkh$\ddot{a}$user
Verlag, Basel,
 1992.



\bibitem {L} Sahu, L. :
 ``Quantum Stochastic Dilation of a Class  of Quantum Dynamical
Semigroups and Quantum Random Walks", Ph. D. Thesis, Indian
Statistical Institute,  2005.



\bibitem {S} Sahu, L. :
Quantum random walks and their convergence, Submitted to ``Infinite
Dimensional Analysis and Quantum Probability and Related Topic"
2005. Available at:
http://xxx.lanl.gov/abs/math.OA/0505438 \\

\bibitem {SiSm} Sinclair, Allan M. and  Smith, Roger R. :
``Hochschild cohomology of von Neumann algebras",
 London Mathematical Society Lecture Note Series, {\bf 203},
 Cambridge University Press, Cambridge, 1995.

\bibitem {KBS} Sinha, K.B. :  Quantum random walks revisited, submitted to
Proc. QP conference, Bedlewo, Poland, June-2004.



\end{thebibliography}
\end{document}